\documentclass[a4paper,12pt]{amsart}
\usepackage{kv2}

\begin{document}

\title{Nonwandering sets of interval skew products}

\author{V.~Kleptsyn}
\address{Victor~Kleptsyn
\newline\hphantom{iii} CNRS, Institut de Recherche Mathematique de Rennes (IRMAR, UMR 6625 CNRS)}
\email{Victor.Kleptsyn@univ-rennes1.fr}
\thanks{V.\ K.\ was supported in part by RFBR project 13-01-00969-a and by a RFBR/CNRS joint project 10-01-93115-CNRS\_a}

\author{{D.~Volk}}
\address{Denis~Volk
\newline\hphantom{iii} University of Rome ``Tor Vergata''
\newline\hphantom{iii} Institute for Information Transmission Problems, Russian Academy of Sciences}
\curraddr{Dipartimento di Matematica, Via della Ricerca Scientifica, 00133 Roma Italy}
\email{volk@mat.uniroma2.it}
\thanks{D.\ V.\  was supported in part by European Advanced Grant MALADY (ERC AdG 246953), and grants RFBR 12-01-31241-mol\_a, RFBR 13-01-00969-a, President's of Russia MK-7567.2013.1}




\begin{abstract}
In this paper we consider a class of skew products over transitive subshifts of finite type with interval fibers. For a natural class of 1-parameter families we prove that for all but countably many parameter values the nonwandering set (in particular, the union of all attractors and repellers) has zero measure. As a consequence, the same holds for a residual subset of the space of skew products.
\end{abstract}

\maketitle

\section{Introduction}  \label{s:intro}

Skew products over hyperbolic dynamics or, equivalently, over subshifts of finite type is quite a standard component of the modern theory of partially hyperbolic dynamical systems. A detailed review of this role of skew products can be found, for instance, in 
the introduction to~\cite{Kleptsyn2014a}. For more references on partial hyperbolicity one can check classical works~\cite{Brin1974}, \cite{HPS1977}, and more modern~\cite{Bonatti2004}.

In~\cite{Kleptsyn2014a} which is the prequel to this paper, we started
a project to understand such skew products in the case when the fibers are intervals of real line, and the fiber maps are orientation preserving diffeomorphisms. The first milestone on this road was to describe the dynamics of so-called \emph{step} skew products: such that the dependence of fiber maps on the base coordinate is piecewise constant.

It is well known that in~$\dim = 1$, that is, circle or interval, generic diffeomorphisms are Morse-Smale: the nonwandering set is a finite collection of periodic points each of whose is hyperbolic. In the case of an orientation-preserving diffeomorphism of a closed interval, this is enough to give the complete description of the dynamics: the nonwandering set is a finite collection of fixed points, each of whose is either attracting or repelling. Attractors and repellers come interchanging: ARAR\dots RA, and all the orbits from any $(A,R)$-interval are monotone moving from the repeller to the attractor.

In~\cite{Kleptsyn2014a} we showed that there exists an open and dense subset of  the set of step skew products where the dynamics is very similar to the \emph{direct} product of the subshift in the base and a generic diffeomorphism of the unit interval. In particular, almost all the nonwandering set consists of finitely many attractors and repellers which are invariant graphs of some continuous maps from the base to the fiber. Attractors and repellers come interchanging, and all the orbits starting between given attractor and repeller drift towards the attractor in the sense of Definition~\ref{d:point_drifts}.

In this sequel we manage to relax the \emph{step} condition and to extend some of the results of~\cite{Kleptsyn2014a} to \emph{all} skew products. Namely, in Theorem~\ref{t:small} we show that for a residual set of such skew products, almost every point with respect to the \emph{standard measure} is either drifting up or down. Here the standard measure is the product of an ergodic invariant measure of the subshift in the base and the Lebesgue measure in the fiber. This implies that the nonwandering set (including all the attractors and repellers) has zero standard measure.

For this, we prove a stronger statement which is Theorem~\ref{t:countable}: for any monotone 1-parameter family, for all but countably many parameter values the standard measure of so-called \emph{anchored} set (which includes the nonwandering set) is zero.

The authors thank the anonymous referee for the careful reading and numerous valuable comments and suggestions which helped to improve the paper.

\section{Notations} \label{s:notations}

Let~$\sigma \colon \Sigma \to \Sigma$ be a bilateral transitive subshift of finite type (a topological Markov chain), $\Sigma \subset \{1, \ldots, N\}^\bbZ$. We endow~$\Sigma$ with the metric defined by the formula
\begin{equation}    \label{e:markov_metric}
d (\bar\om, \om) = \begin{cases}
2^{-\min\{|n| \, : \, \bar\om_n\neq \om_n \}}, &  \bar\om \neq \om, \\
0, &   \bar\om = \om,
\end{cases}
\qquad \bar\om,\om \in \Sigma.
\end{equation}
Fix any $\sigma$-invariant ergodic Markov measure~$\mu$ on $\Sigma$.

Let $I \subset \bbR$ be a unit interval.
In this paper, we will be considering the following class~$\mcS$ of skew products $F \colon \Sigma \times I \to \Sigma \times I$:
\begin{enumerate}
\item\label{e:skpr-class-1} $F \colon (\om, x) \mapsto (\sigma\om, f_\om(x))$;
\item\label{e:skpr-class-2} the maps~$f_\om \colon I \to f_\om(I) \subset I$ are orientation preserving $C^1$ diffeomorphisms;
\item\label{e:skpr-class-3} the map~$f_{(\cdot)} \colon \Sigma \to \Diff^1(I)$ is continuous.
\end{enumerate}
We equip~$\mcS$ with a metric as follows:
\begin{equation} \label{e:metric}
\dist(F, \tilde F) = \max_{\om} \dist_{C^1} (f_\om^{\pm 1}, \tilde f_\om^{\pm 1}).
\end{equation} 
Note that for an open and dense subset of~$\mcS$, the interval~$I$ is mapped \emph{strictly} inside itself. 

The \emph{standard measure} $\mathbf{m}$ on $\Sigma \times I$ is the product of measure~$\mu$ in the base and the Lebesgue measure in the fiber.
Finally, let~$I_\omega \subset \Sigma \times I$ be the fiber over~$\omega \in \Sigma$, and for any subset~$K \subset \Sigma \times I$ denote $K_\omega = K \cap I_\omega$.

\section{Drifting and anchored regions}    \label{s:drifting-regions}

Following~\cite{Kleptsyn2014a}, we define drifting graphs and points.
Let $\varphi_i \colon \Sigma \to I$ be two arbitrary functions, $\Gamma_i$ be their graphs, $i = 1,2$.

\begin{Def} \label{d:compare}
We write $\varphi_1 < \varphi_2$ whenever for any $\om \in \Sigma$
$$
\quad \varphi_1(\om) < \varphi_2 (\om).
$$
We also write~$\Gamma_1 < \Gamma_2$ in this case.
\end{Def}

\begin{Def} \label{d:band}
A \emph{horizontal band} is the set of the points between~$\Gamma_1$ and~$\Gamma_2$, $\Gamma_1 < \Gamma_2$.
\end{Def}

A general property of any skew product with an invertible map in the base is that the image~$F(\Gamma)$ of any graph~$\Gamma$ is also a graph of some function.

\begin{Def} \label{d:drift_up}
We say that a graph~$\Gamma$ \emph{drifts up (down)} if $F(\Gamma) > \Gamma$ (respectively, $F(\Gamma) < \Gamma$).
\end{Def}

Recall that we assume the fiber maps~$f_\om$ to be monotone increasing. Thus $F(\Gamma) > \Gamma$ implies~$F^{n+1}(\Gamma) > F^n(\Gamma)$ for any~$n \in \bbZ$, and $F(\Gamma) < \Gamma$ implies~$F^{n+1}(\Gamma) < F^n(\Gamma)$ for any~$n \in \bbZ$, provided the backward iterates are well-defined. The band between $F^n(\Gamma)$ and $F^{n+1}(\Gamma)$ is a fundamental domain for $F$ restricted to a certain horizontal band, see Section~\ref{s:geom} for more details.

\begin{Def} \label{d:point_drifts}
A point $p=(\omega,x)\in \Sigma\times I$ \emph{drifts up}
if there exists a continuous map~$\gamma \colon \Sigma\to I$ with the graph~$\Gamma$ such that
\begin{itemize}
\item $\Gamma$ drifts up under~$F$;
\item the point $p$ is between~$\Gamma$ and its image:
$$
\gamma(\omega)<x<f_{\sigma^{-1}\omega}(\gamma(\sigma^{-1}\omega)).
$$
\end{itemize}
Denote the set of the points drifting up in~$F$ by $\Up(F)$. In the same way we define the set of points \emph{drifting down} which we denote by~$\Down(F)$. Finally, denote by~$\Indiff(F)=(\Sigma \times I)\setminus (\Up(F)\cup \Down(F))$ the \emph{anchored set} which is the complement to the points drifting up or down.
\end{Def}

\begin{Rem}	\label{r:omega}
For the nonwandering set $\Omega(F)$ we always have~$\Omega(F) \subset \Indiff(F)$ because the sets $\Up(F)$ and $\Down(F)$ are open and disjoint from~$\Omega(F)$.
\end{Rem}

\begin{Prop}    \label{p:up-down-ni}
For any~$F$, we have $\Up(F) \cap \Down(F) = \emptyset$.
\end{Prop}

\begin{proof}
Assume that $\Up(F)\cap \Down(F) \ne \emptyset$ for some skew product~$F$. Because the set $\Up(F)\cap \Down(F)$ is open and the periodic points are dense in~$\Sigma$, we can pick~$p = (\omega, x) \in \Up(F)\cap \Down(F)$ such that~$\omega$ is periodic. Let~$n$ be its minimal period. Then~$F^n(p) = (\bar x, \bar\omega)$ belongs to the same fiber as~$p$. Note that if $\bar x > x$, then~$p$ cannot be drifting down, and if $\bar x < x$, then $p$ cannot be drifting up. Thus $p \notin \Up(F)\cap \Down(F)$, and the Proposition is proven.
\end{proof}

Now we introduce a partial ordering on the set~$\mcS$:

\begin{Def} \label{def:order}
$F \prec \tilde F$ whenever for any $\om,x$ we have $f_\om(x) < \tilde f_\om (x)$.
\end{Def}

\begin{Prop} \label{p:compare}
If $F \prec \tilde F$, then
$$
\Up(F) \subset \Up(\tilde F) \text{ and } \Down(F) \supset \Down(\tilde F).
$$
\end{Prop}

\begin{proof}
  Indeed, for any $p \in \Up(F)$ we can take~$\Gamma$ that satisfies Definition~\ref{d:point_drifts}. Because~$F \prec \tilde F$, we have $F(\Gamma) < \tilde F (\Gamma)$. Thus the same~$\Gamma$ is also valid to show~$p \in \Up(\tilde F)$.

  The second inclusion is proved in the same way.
\end{proof}

\begin{Def}
A family $F_{\tau}$, $\tau \in (\tau_1, \tau_2)$, of skew products is \emph{monotone increasing} if for any $\tau_1 < \tau_2$ the skew products~$F_{\tau_1}$ and $F_{\tau_2}$ are comparable and $F_{\tau_1} \prec F_{\tau_2}$.
\end{Def}

Definition of monotone decreasing family is analogous. A family is monotone of it is either monotone increasing or decreasing.


We say a family~$F_\tau$ is continuous if the map~$\tau\mapsto F_\tau$ is continuous with respect to the metric~\eqref{e:metric}.

\section{Main results}	\label{s:main-results}

\begin{Thm}\label{t:countable}
Let $F_{\tau}$ be a monotone family, continuous in $\tau$. Then for every~$\tau$, except for at most countable set of them, we have
\begin{equation}\label{eq:UpDownCompl}
\mathbf{m}( \Indiff(F_{\tau}))=0,
\end{equation}
where~$\mathbf{m}$ is the standard measure.
In particular, the standard measure of the nonwandering set~$\Omega(F_{\tau}) \subset \Indiff(F_{\tau})$ is zero.

Moreover, for any $\varepsilon>0$ the set $\{\tau\mid\mathbf{m}( \Indiff(F_{\tau})) \ge \eps\}$ is finite.
\end{Thm}

\begin{Def}
We say that a subset of the space of skew products is \emph{small} if it is closed, nowhere dense, and any continuous monotone family intersects it at a finite number of points.
\end{Def}

\begin{Thm}\label{t:small}
The set $\mathcal{K} = \{ F \mid \mathbf{m}(\Indiff(F))>0 \}$ is a subset of a countable union of small sets. In particular, $\mathcal{K}$ is meager.
\end{Thm}

\begin{Rem}	\label{r:fubini}
Combined with Remark~\ref{r:omega}, this implies that a generic skew product from~$\mathcal{S}$ has the nonwandering set~$\Omega$ of zero standard measure. By Fubini's Theorem, for almost every~$\omega \in \Sigma$ the restriction~$\Omega_\omega$ of~$\Omega$ to the fiber over~$\omega$ has zero Lebesgue measure. 
\end{Rem}

\begin{Rem}	\label{r:genericity}
The word ``generic'' is essential here, because some trivial examples of skew products \emph{do} belong to the ``bad'' set~$\mathcal{K}$. In particular, such is the skew product with all identity maps in the fibers. Its nonwandering set is the whole phase space.
\end{Rem}

\begin{Rem}	\label{r:all-markov}
We would also like to emphasize here that~$\mathcal{K}$ does \emph{depend}, in general, on the initial choice of Markov measure in the base. We believe that it is still possible to find a uniform meager set~$\tilde{\mathcal{K}}$ 
which is \emph{simultaneously} valid for the whole continuum of Markov measures in the base, but the argument is beyond the scope of this paper.
\end{Rem}

In the spirit of the Large Deviations Lemma~\cite[Lemma~6]{Ilyashenko2008} and the Special Ergodic Theorem~\cite[Theorem~6]{Ilyashenko2008} which give an estimate of the Hausdorff dimension of ``bad'' sets in some partially hyperbolic systems, we conjecture the following generalization of Theorem~\ref{t:countable}.

\begin{Con} \label{con:countable_hausdorff}
Let $F_{\tau}$ be a monotone family, continuous in $\tau$. Then for every~$\tau$, except for at most countable set of them, the Hausdorff dimension of~$\Indiff(F_{\tau})$ is less than the full dimension of the phase space.
\end{Con}

Finally, we would like to state a very general conjecture. It is based on the Baxendale-like results which imply strong bunching of orbits for skew products with 1-dimensional fibers, see~\cite{Kleptsyn2004}, \cite{Kleptsyn2014a} for applications or~\cite{Baxendale1989} for the general theorem. Additionally, for interval fibers the monotonicity of the fiber maps provides very strong vertical ordering of all the invariant sets/measures, as we will see in use in this paper.

So, fix any tuple~$(X, h, \nu)$, where~$X$ is a compact topological space, $h \colon X \to X$ is a homeomorphism with dense periodic orbits (see Proposition~\ref{p:up-down-ni}), and~$\nu$ is an ergodic invariant measure with full support. Consider the class~$\mathcal{G}$ of skew products with the base~$(X, h, \nu)$, and the fiber maps~$f_\omega$, $\omega \in X$, satisfying the same conditions~(\ref{e:skpr-class-2}), (\ref{e:skpr-class-3}) from Section~\ref{s:notations} as before. Analogously, the standard measure is the product of the ergodic measure~$\nu$ in the base and the Lebesgue measure in the fibers.

\begin{Con} \label{con:general}
Let $G_{\tau} \in \mathcal{G}$ be a monotone family of skew products, continuous in $\tau$.  Then for every~$\tau$, except for at most countable set of them, the anchored set has zero standard measure.
\end{Con}

Back to the ground, in the following sections we first discuss the geometrical implications of Theorems~\ref{t:countable} and~\ref{t:small}, and then prove them in Sections~\ref{s:t:countable} and~\ref{s:t:small}.

\section{Geometry of drifting and anchored regions}	\label{s:geom}

In this Section, we discuss geometrical implications of Theorem~\ref{t:countable}. They are rather straightforward but we find them worth mentioning, especially in comparison with our results~\cite{Kleptsyn2014a}.

Take any~$p \in \Up(F)$ and the corresponding horizontal band~$\mcB \ni p$ between some continuous graphs~$\Gamma$ and $F(\Gamma)$, $\Gamma < F(\Gamma)$. The bands~$F^n(\mcB)$, $n \in \bbZ$, are pairwise disjoint and naturally ordered:
$$
F^i(\mcB) < F^j (\mcB) \quad \text{for} \quad i <j.
$$
We call~$\mcD = \mcD(p) = \cup_{n \in \bbZ} F^n (\mcB) \subset \Up(F)$ the \emph{maximal component of $p$} in $\Up(F)$. $\mcD$ is an invariant open band bounded from below and above by the monotone pointwise limits
$$
D_- = \lim_{n \to -\infty} F^n (\Gamma) \quad \text{and} \quad D_+ = \lim_{n \to +\infty} F^n (\Gamma),
$$
respectively. By standard calculus, $D_-$ is the graph of an upper semicontinuous function~$d_-(\omega)$, and $D_+$ is the graph of a lower semicontinuous function~$d_+(\omega)$. Note that $d_-(\omega)$ and $d_+(\omega)$ are $F$-invariant sections: $f_\omega(d_{\pm}(\omega)) = d_\pm (\sigma \omega)$. By construction, $\mcD_\omega$ is an open interval $(d_-(\omega), d_+(\omega)$ for any~$\omega \in \Sigma$. So for any~$\omega$, $\Up_\omega$ is a collection of disjoint intervals, thus a countable one.


Obviously, the structure of~$\Down(F)$ is analogous.
However, the geometry of~$\Indiff(F)$ may be more complex. Because the set~$\Indiff$ is the complement to~$\Up \cup \Down$, we know that~$\Indiff_\omega$ is a closed set for any~$\omega \in \Sigma$. Also, $F(\Indiff_\omega) = \Indiff_{\sigma\omega}$. Recall that in~\cite{Kleptsyn2014a} we showed that for an open and dense set of \emph{step} skew products $\Indiff_\omega$ is a disjoint union of finitely many singleton points or closed intervals. In the setting of present paper, however, $\Indiff_\omega$ may also contain Cantor-like subsets.

\section{Proof on Theorem~\ref{t:countable}}	\label{s:t:countable}

We give the proof for monotone increasing~$F_\tau$. The case of monotone decreasing family is handled by replacing~$\tau \mapsto -\tau$.

\begin{Def}
A skew product~$F \colon \Sigma \times M \to \Sigma \times M$ is \emph{multistep} if the fiber maps~$f_\om$ depend only on finitely many positions in~$\om$.
\end{Def}

In the same way every continuous function can be approximated by piecewise constant functions in $\sup$-norm,  we can $C^0 (\Sigma, C^1(I))$-approximate skew products from~$\mcS$ by multistep skew products. These approximations can be chosen to be generic in the sense of \cite{Kleptsyn2014a}. 
It immediately follows from~\cite[Theorem~2.15]{Kleptsyn2014a} that for any generic multistep skew product $G$ we have $\mathbf{m}(\Indiff(G)) = 0$.

\begin{Prop}    \label{p:comparison}
Let $F_1 \prec F_2$. Then $\mathbf{m}(\Down(F_1) \cup \Up(F_2)) = 1$.
\end{Prop}

\begin{proof}
Take a generic multistep skew product~$G$ such that $F_1 \prec G \prec F_2$. Because of the above remark,
$$
\mathbf{m}(\Down(G) \cup \Up(G)) = 1.
$$
But $\Down(F_1) \supset \Down(G)$ and $\Up(F_2) \supset \Up(G)$. Thus~$\Down(G) \cup \Up(G) \subset \Down(F_1) \cup \Up(F_2)$ which proves the Proposition.
\end{proof}

Now fix any monotone increasing family~$F_{\tau}$.
Denote for brevity $\Up_{\tau}:=\Up(F_{\tau})$, $\Down_{\tau}:=\Down(F_{\tau})$.
Proposition~\ref{p:compare} implies that for any $\tau_1 < \tau_2$ we have $\Up_{\tau_1} \subset \Up_{\tau_2}$ and $\Down_{\tau_1} \supset \Down_{\tau_2}$.
Also denote
$$
\Up_{\tau+}:=\bigcap_{\delta>0} \Up_{\tau+\delta}, \quad \Down_{\tau-}:=\bigcap_{\delta>0} \Down_{\tau-\delta}.
$$
Obviously, $\Up_\tau \subset \Up_{\tau+}$ and $\Down_\tau \subset \Down_{\tau-}$.

\begin{Prop}\label{p:measure}
For any~$\tau$,
\begin{enumerate}
\item\label{i:plus} $\Up_{\tau+} \cap \Down_{\tau} = \emptyset$.
\item\label{i:measure} $\mathbf{m}(\Up_{\tau+}\cup \Down_{\tau})=1$.
\end{enumerate}
\end{Prop}
\begin{proof}
Because~$F_\tau$ is continuous in~$\tau$, for any point $p$ the set of parameters $\{\tau \mid p\in \Down_{\tau} \}$ is open. Thus for any $p\in \Down_{\tau}$ for any small enough $\delta>0$ we have $p\in \Down_{\tau+\delta}$. By Proposition~\ref{p:up-down-ni}, this implies $p\notin \Up_{\tau+\delta}$. Taking the intersection over all~$\delta > 0$, we have that $p$ does not belong to~$\Up_{\tau+}$. Thus $\Up_{\tau+} \cap \Down_{\tau} = \emptyset$.

Let us now prove~\ref{i:measure}.
%
By Proposition~\ref{p:comparison}, for any $\delta>0$
$$
\mathbf{m}\left(\Down_{\tau}\bigcup \Up_{\tau+\delta}\right) = 1.
$$
By Proposition~\ref{p:compare}, the sets~$\Up_{\tau+\delta}$ are monotone increasing in~$\delta$. Take the intersection over all~$\delta > 0$ to get
$$
\mathbf{m}\left(\bigcap_{\delta>0}\left(\Down_{\tau}\cup \Up_{\tau+\delta}\right)\right) = 1.
$$
Now factor out the term~$\Down_{\tau}$ to get the required $\mathbf{m}(\Down_{\tau}\cup \Up_{\tau+})=1$.
\end{proof}

Now we are ready to complete the
\begin{proof}[Proof of Theorem~\ref{t:countable}]
Because $\Down_\tau \cap \Up_{\tau+} = \emptyset$ and $\mathbf{m}(\Down_{\tau}\cup\Up_{\tau+}) = 1$,
we have
$$
\mathbf{m}(\Indiff(F_{\tau})) = \mathbf{m}((\Sigma\times [0,1])\setminus (\Down_{\tau} \sqcup \Up_{\tau}))= \mathbf{m}((\Down_{\tau} \sqcup \Up_{\tau+})\setminus (\Down_{\tau} \sqcup \Up_{\tau})) = \mathbf{m}(\Up_{\tau+}\setminus \Up_{\tau}).
$$
This can be rewritten as
$$
\mathbf{m}(\Up_{\tau+}\setminus \Up_{\tau}) = \mathbf{m}\left(\bigcap_{\delta>0} \Up_{\tau+\delta}\right) - \mathbf{m}(\Up_{\tau}) = \lim_{\delta\to +0} \mathbf{m}(\Up_{\tau+\delta}) - \mathbf{m}(\Up_{\tau}).
$$
So $\mathbf{m}(\Indiff(F_{\tau}))$ equals to the value of the gap of the monotone increasing function~$\mu(t)=\mathbf{m}(\Up_t)$ at the point $t = \tau$. Any  monotone function has at most countable number of gaps. Moreover, if the function is bounded, then for any fixed~$\eps > 0$ only finitely many of the gaps can be bigger than~$\eps$. Theorem~\ref{t:countable} is proven.
\end{proof}

Finally, we remark that 
all continuous monotone increasing families $F_\tau, G_\rho$ such that $F_0 = G_0$ are equivalent in the following sense: for any $\tau > 0$ there exists $\rho > 0$ such that 
$$
F_0 \prec G_\rho \prec F_\tau.
$$
%
Indeed, one can take $0 < \eps = \inf_{\om,x}(F_\tau(\om,x) - F_0(\om,x))$, and take $\rho > 0$ such that $0 < G_\rho - G_0 < \eps$. Then $F_0 = G_0 \prec G_\rho \prec F_\tau$.
%
Because this argument is symmetric with respect to switching $F$ and $G$,  the sets $\Up_{0+}, \Down_{0-}$ depend only on~$F_0$ but not on the choice of a continuous monotone increasing family passing through~$F_0$. Thus we could denote them just by $\Up^+ (F_0)$,  $\Down^- (F_0)$.

\section{Proof of Theorem~\ref{t:small}}	\label{s:t:small}
Obviously, $\mathcal{K} = \cup_{n \in \bbN} \mathcal{K}_{1/n}$, where
$$
\mathcal{K}_\eps := \{ F \mid \mathbf{m}(\Indiff(F)) \ge \eps \}.
$$
Let us show that for any~$\eps > 0$ the set $\mathcal{K}_\eps$ is small. First of all, Theorem~\ref{t:countable} implies that any monotone family intersects~$\mathcal{K}_\eps$ at a finite number of points.

\begin{Prop}
For any~$\eps > 0$, the set $\mathcal{K}_\eps$ is closed.
\end{Prop}

\begin{proof}
Note that the set of pairs $\mathcal{S}_{\Up}=\{(p,F) \mid p\in \Up(F) \}$ is an open subset of the Cartesian product of $\Sigma \times I$ and the space of dynamical systems on it. The same is true for $\mathcal{S}_{\Down}$. Thus the set
$$
\mathcal{S}_{\Indiff}:= \{(p,F)\mid p\in \Indiff(F) \} = \{(p,F)\mid p\notin \Up(F)\cup \Down(F) \}
$$
is closed.

Now take any sequence of systems~$F_{(n)} \in \mathcal{K}_\eps$ which converges to~$F$. Assume that~$F \notin \mathcal{K}_\eps$ which means $\mathbf{m}(\Indiff(F)) < \eps$. Take an open cover~$U$ of the set~$\Indiff(F)$ such that~$\mathbf{m}(U) < \eps$. Because $F_{(n)} \in \mathcal{K}_\eps$, we have $\mathbf{m}(\Indiff(F_{(n)})) \ge \eps$ for all $n \in \bbN$. Thus~$\Indiff(F_{(n)}) \not\subset U$. Then for any $n \in \bbN$ there exists a point~$p_n \in (\Sigma \times I) \setminus U$, $p_n \in \Indiff(F_{(n)})$.

Because the set $(\Sigma \times I) \setminus U$ is compact, we can extract from~$p_n$ a converging subsequence~$p_{n_m}$. Then the subsequence $(p_{n_m}, F_{(n_m)}) \in \mathcal S_\Indiff$ converges to some $(\tilde p, F)$, and $\tilde p \notin U$. But because the set $\mathcal S_\Indiff$ is closed, we must have $p \in \Indiff(F)$. The contradiction with $\Indiff(F) \subset U$ proves the Proposition.
\end{proof}

Because any monotone family intersects~$\mathcal{K}_{\varepsilon}$ at a finite number of points, the set $\mathcal{K}_{\varepsilon}$ has empty interior. Thus it is nowhere dense. Theorem~\ref{t:small} is proven.

\bibliographystyle{plain}
\bibliography{dynsys-utf8}

\end{document}